\documentclass[12pt]{article}
\usepackage{times,amsmath,amssymb,theorem,xspace}
\usepackage[dvips]{graphicx}
\usepackage[hscale=0.73,vscale=0.79]{geometry}

\theoremheaderfont{\scshape}
\theorembodyfont{\normalfont\slshape}
\theoremstyle{plain}
\newtheorem{theorem}{Theorem}
\newtheorem{proposition}{Proposition}
\newtheorem{lemma}{Lemma}

\theorembodyfont{\normalfont}

\newenvironment{proof}{\begin{trivlist}\item{}\normalfont\textit{Proof. }}{\hfill$\square$\end{trivlist}}
\newenvironment{proofof}[1]{\begin{trivlist}\item{}\normalfont\textit{Proof of #1.}}{\hfill$\square$\end{trivlist}}

\newdimen\proofrulebreadth \proofrulebreadth=.05em
\newdimen\proofdotseparation \proofdotseparation=1.25ex
\newdimen\proofrulebaseline \proofrulebaseline=2ex
\newcount\proofdotnumber \proofdotnumber=3
\let\then\relax
\def\hfi{\hskip0pt plus.0001fil}
\mathchardef\squigto="3A3B
\newif\ifinsideprooftree\insideprooftreefalse
\newif\ifonleftofproofrule\onleftofproofrulefalse
\newif\ifproofdots\proofdotsfalse
\newif\ifdoubleproof\doubleprooffalse
\let\wereinproofbit\relax
\newdimen\shortenproofleft
\newdimen\shortenproofright
\newdimen\proofbelowshift
\newbox\proofabove
\newbox\proofbelow
\newbox\proofrulename
\def\shiftproofbelow{\let\next\relax\afterassignment\setshiftproofbelow\dimen0 }
\def\shiftproofbelowneg{\def\next{\multiply\dimen0 by-1 }%
\afterassignment\setshiftproofbelow\dimen0 }
\def\setshiftproofbelow{\next\proofbelowshift=\dimen0 }
\def\setproofrulebreadth{\proofrulebreadth}

\def\prooftree{
\ifnum  \lastpenalty=1
\then   \unpenalty
\else   \onleftofproofrulefalse
\fi
\ifonleftofproofrule
\else   \ifinsideprooftree
        \then   \hskip.5em plus1fil
        \fi
\fi
\bgroup
\setbox\proofbelow=\hbox{}\setbox\proofrulename=\hbox{}%
\let\justifies\proofover\let\leadsto\proofoverdots\let\Justifies\proofoverdbl
\let\using\proofusing\let\[\prooftree
\ifinsideprooftree\let\]\endprooftree\fi
\proofdotsfalse\doubleprooffalse
\let\thickness\setproofrulebreadth
\let\shiftright\shiftproofbelow \let\shift\shiftproofbelow
\let\shiftleft\shiftproofbelowneg
\let\ifwasinsideprooftree\ifinsideprooftree
\insideprooftreetrue
\setbox\proofabove=\hbox\bgroup$\displaystyle 
\let\wereinproofbit\prooftree
\shortenproofleft=0pt \shortenproofright=0pt \proofbelowshift=0pt
\onleftofproofruletrue\penalty1
}

\def\eproofbit{
\ifx    \wereinproofbit\prooftree
\then   \ifcase \lastpenalty
        \then   \shortenproofright=0pt  
        \or     \unpenalty\hfil         
        \or     \unpenalty\unskip       
        \else   \shortenproofright=0pt  
        \fi
\fi
\global\dimen0=\shortenproofleft
\global\dimen1=\shortenproofright
\global\dimen2=\proofrulebreadth
\global\dimen3=\proofbelowshift
\global\dimen4=\proofdotseparation
\global\count255=\proofdotnumber
$\egroup  
\shortenproofleft=\dimen0
\shortenproofright=\dimen1
\proofrulebreadth=\dimen2
\proofbelowshift=\dimen3
\proofdotseparation=\dimen4
\proofdotnumber=\count255
}

\def\proofover{
\eproofbit 
\setbox\proofbelow=\hbox\bgroup 
\let\wereinproofbit\proofover
$\displaystyle
}%
\def\proofoverdbl{
\eproofbit 
\doubleprooftrue
\setbox\proofbelow=\hbox\bgroup 
\let\wereinproofbit\proofoverdbl
$\displaystyle
}%
\def\proofoverdots{
\eproofbit 
\proofdotstrue
\setbox\proofbelow=\hbox\bgroup 
\let\wereinproofbit\proofoverdots
$\displaystyle
}%
\def\proofusing{
\eproofbit 
\setbox\proofrulename=\hbox\bgroup 
\let\wereinproofbit\proofusing
\kern0.3em$
}

\def\endprooftree{
\eproofbit 
  \dimen5 =0pt
\dimen0=\wd\proofabove \advance\dimen0-\shortenproofleft
\advance\dimen0-\shortenproofright
\dimen1=.5\dimen0 \advance\dimen1-.5\wd\proofbelow
\dimen4=\dimen1
\advance\dimen1\proofbelowshift \advance\dimen4-\proofbelowshift
\ifdim  \dimen1<0pt
\then   \advance\shortenproofleft\dimen1
        \advance\dimen0-\dimen1
        \dimen1=0pt
        \ifdim  \shortenproofleft<0pt
        \then   \setbox\proofabove=\hbox{%
                        \kern-\shortenproofleft\unhbox\proofabove}%
                \shortenproofleft=0pt
        \fi
\fi
\ifdim  \dimen4<0pt
\then   \advance\shortenproofright\dimen4
        \advance\dimen0-\dimen4
        \dimen4=0pt
\fi
\ifdim  \shortenproofright<\wd\proofrulename
\then   \shortenproofright=\wd\proofrulename
\fi
\dimen2=\shortenproofleft \advance\dimen2 by\dimen1
\dimen3=\shortenproofright\advance\dimen3 by\dimen4
\ifproofdots
\then
        \dimen6=\shortenproofleft \advance\dimen6 .5\dimen0
        \setbox1=\vbox to\proofdotseparation{\vss\hbox{$\cdot$}\vss}%
        \setbox0=\hbox{%
                \advance\dimen6-.5\wd1
                \kern\dimen6
                $\vcenter to\proofdotnumber\proofdotseparation
                        {\leaders\box1\vfill}$%
                \unhbox\proofrulename}%
\else   \dimen6=\fontdimen22\the\textfont2 
        \dimen7=\dimen6
        \advance\dimen6by.5\proofrulebreadth
        \advance\dimen7by-.5\proofrulebreadth
        \setbox0=\hbox{%
                \kern\shortenproofleft
                \ifdoubleproof
                \then   \hbox to\dimen0{%
                        $\mathsurround0pt\mathord=\mkern-6mu%
                        \cleaders\hbox{$\mkern-2mu=\mkern-2mu$}\hfill
                        \mkern-6mu\mathord=$}%
                \else   \vrule height\dimen6 depth-\dimen7 width\dimen0
                \fi
                \unhbox\proofrulename}%
        \ht0=\dimen6 \dp0=-\dimen7
\fi
\let\doll\relax
\ifwasinsideprooftree
\then   \let\VBOX\vbox
\else   \ifmmode\else$\let\doll=$\fi
        \let\VBOX\vcenter
\fi
\VBOX   {\baselineskip\proofrulebaseline \lineskip.2ex
        \expandafter\lineskiplimit\ifproofdots0ex\else-0.6ex\fi
        \hbox   spread\dimen5   {\hfi\unhbox\proofabove\hfi}%
        \hbox{\box0}%
        \hbox   {\kern\dimen2 \box\proofbelow}}\doll%
\global\dimen2=\dimen2
\global\dimen3=\dimen3
\egroup 
\ifonleftofproofrule
\then   \shortenproofleft=\dimen2
\fi
\shortenproofright=\dimen3
\onleftofproofrulefalse
\ifinsideprooftree
\then   \hskip.5em plus 1fil \penalty2
\fi
}

\makeatletter
\def\section{\@startsection{section}{1}{0pt}{-3.25ex plus -1ex minus 
-.2ex}{1.5ex plus .2ex minus .3ex}{\normalfont\large\bf}}
\def\subsection{\@startsection {subsection}{2}{0pt}{-2ex plus -1ex minus 
   -.2ex}{1.5ex plus .2ex minus .3ex}{\normalfont\normalsize\bf}}
\makeatother

\newlength{\tw}
\newenvironment{labelledbox}[1]{\begin{flushleft}\vspace*{-2.6ex}\rule{5mm}{0.4pt}%
\hspace*{0.3mm}\raisebox{-3mm}{\rule{0mm}{6mm}}
\raisebox{-0.8mm}{\bf #1}{\rule{0mm}{6mm}}
\hrulefill\vspace*{-3.4mm}
\begin{tabular}{@{}|c|@{}}\begin{minipage}[b]{\tw}\vspace*{5.4mm}\begingroup}
{\protect\vspace*{-1mm}\endgroup\end{minipage}\\
 \hline \end{tabular}\end{flushleft}}

\newcommand{\defn}[1]{{\textit{\textbf{#1}}}}
\newcommand{\myitem}[1]{\item[\textnormal{(#1)}]}

\newcommand{\ie}{\emph{i.e.}}
\newcommand{\etc}{\emph{etc.}}
\newcommand{\eg}{\emph{e.g.}}
\newcommand{\Cf}{\emph{Cf.}}
\newcommand{\cf}{\emph{cf.}}

\newcommand{\with}{\&}
\newcommand{\plus}{\oplus}
\newcommand{\tensor}{\otimes}

\newlength{\parrdp}\newlength{\parrht}
\newcommand{\parr}{\raisebox{-\parrdp}{\raisebox{\parrht}{\rotatebox{180}{$\&$}}}}
\newlength{\smallparrdp}\newlength{\smallparrht}
\newcommand{\smallparr}
{\mkern1mu\raisebox{-.1\smallparrdp}{\raisebox{\smallparrht}{\rotatebox{180}{\small$\&$}}}\mkern1mu}
\newlength{\footnoteparrdp}\newlength{\footnoteparrht}
\newcommand{\footnoteparr}
{\mkern1mu\raisebox{-\footnoteparrdp}{\raisebox{\footnoteparrht}{\rotatebox{180}{\footnotesize$\&$}}}\mkern1mu}
\newlength{\scriptparrdp}\newlength{\scriptparrht}
\newlength{\tinyparrdp}\newlength{\tinyparrht}
\newcommand{\weakening}{\mathsf{W}}
\newcommand{\contraction}{\mathsf{C}}
\newcommand{\rulelabel}[1]{\mathsf{#1}}
\newcommand{\axlabel}{}
\newcommand{\contractionlabel}{\;\contraction}
\newcommand{\weakeninglabel}{\;\weakening}
\newcommand{\parlabel}{\parr}

\newcommand{\andlabel}{\;\wedge}
\newcommand{\withlabel}{\;\with}
\newcommand{\tensorlabel}{\;\tensor}

\newcommand{\plusilabel}{\;\plus_i}
\newcommand{\plusonelabel}{\;\plus_1}
\newcommand{\plustwolabel}{\;\plus_2}

\newcommand{\cutlabel}{\rulelabel{\mkern1mu cut}}
\newcommand{\supcutlabel}{\rulelabel{\mkern1mu cut}_\with}
\newcommand{\splitcutlabel}{\rulelabel{\mkern1mu cut}_\tensor}

\newcommand{\mixlabel}{\rulelabel{\mkern1mu mix}}
\newcommand{\supmixlabel}{\rulelabel{\mkern1mu mix}_\with}
\newcommand{\splitmixlabel}{\rulelabel{\mkern1mu mix}_\tensor}
\newcommand{\dual}[1]{\overline #1}
\newcommand{\LK}{\textnormal{\textbf{LK}}\xspace}
\newcommand{\core}{\textnormal{\textbf{M}}\xspace}
\newcommand{\corep}{\textnormal{\textbf{Mp}}\xspace}
\newcommand{\gsonep}{\textnormal{\textbf{GS1p}}\xspace}
\newcommand{\gsthree}{\textnormal{\textbf{GS3}}\xspace}
\newcommand{\gsthreep}{\textnormal{\textbf{GS3p}}\xspace}
\newcommand{\gsfivep}{\textnormal{\textbf{GS5p}}\xspace}
\newcommand{\h}{\core}
\newcommand{\Cp}{\corep}
\newcommand{\Cpminus}{\textnormal{\textbf{Mp}$^-$}\xspace}
\newcommand{\positive}{\textnormal{\textbf{Pp}}\xspace}
\newcommand{\negative}{\textnormal{\textbf{Np}}\xspace}
\newcommand{\subforest}{\subseteq}
\newcommand{\contains}{\sqsupseteq}

\title{%
\vspace{-3ex}
{\Large\bf A minimal classical sequent calculus free of structural rules\vspace*{-.5ex}%
}}%
\author{\normalsize
   \sc Dominic Hughes \\ \small Stanford University\thanks{Visiting
   Scholar, Concurrency Group, Computer Science Department, Stanford
   University.  I gratefully acknowledge my host, Vaughan Pratt.}  }
\date{
\small April 22, 2009}
\begin{document}
\maketitle
\thispagestyle{empty}
\vspace*{-5ex}\begin{quotation}\small
Gentzen's classical sequent calculus $\LK$ has explicit
structural rules for contraction and weakening.  They can be absorbed
(in a right-sided formulation) by replacing the axiom ${P,\mkern-2mu
\neg P}$ by ${\mkern1mu\Gamma\!,P,\mkern-2mu \neg P\mkern 1mu}$ for
any context $\Gamma$, and replacing the original disjunction rule with
$\mkern1mu \Gamma\!,A,B\mkern1mu$ implies\/ $
\mkern1mu\Gamma\!,\mkern1mu A\vee B\mkern1mu$.

This paper presents a classical sequent calculus which is also free of
contraction and weakening, but more symmetrically: both contraction
and weakening are absorbed into conjunction, leaving the axiom rule
intact.  It uses a blended conjunction rule, combining the standard
context-sharing and context-splitting rules: $\mkern1mu\Gamma\!,\Delta,A\mkern1mu$
and $\mkern1mu\Gamma\!,\Sigma,B\mkern1mu$ implies\/ $
\mkern1mu\Gamma\!,\Delta,\Sigma,\mkern1mu A\wedge B\mkern1mu$.  
We refer to this system $\core$ as \emph{minimal sequent calculus}.

We prove a \emph{minimality theorem} for the propositional fragment
$\corep$: any propositional sequent calculus $S$ (within a standard
class of right-sided calculi) is complete \emph{if and only if} $S$
contains $\corep$ (that is, each rule of $\corep$ is derivable in
$S$).  Thus one can view $\core$ as a minimal complete core of
Gentzen's $\LK$.
\end{quotation}

\section{Introduction}

The following Gentzen-Sch\"utte-Tait \cite{Gen35,Sch50,Tai68} system,
denoted \gsonep in \cite{TS96}, is a standard right-sided formulation
of the propositional fragment of Gentzen's classical sequent calculus
$\LK$:
\begin{labelledbox}{System \gsonep}{\label{gsonep}\begin{center}
\vspace{1ex}
\begin{prooftree}\thickness=.08em
\rule{0pt}{1.3ex}
\justifies \;P,\neg P\; \using \axlabel
\end{prooftree}
\hspace{8ex}
\begin{prooftree}\thickness=.08em
\Gamma, A\hspace{8ex}\Gamma, B
\justifies \;\Gamma, A\wedge B\; \using \withlabel
\end{prooftree}
\hspace{8ex}
\begin{prooftree}\thickness=.08em
\Gamma, A_i
\justifies \;\Gamma, A_1\vee A_2\; \using \plusilabel
\end{prooftree}

\vspace{4ex}
\begin{prooftree}\thickness=.08em
\Gamma
\justifies \;\Gamma, A \using \weakeninglabel
\end{prooftree}
\hspace{8ex}
\begin{prooftree}\thickness=.08em
\Gamma, A,A
\justifies \;\Gamma, A\; \using \contractionlabel
\end{prooftree}
\vspace{1ex}
\end{center}}\end{labelledbox}
Here $P$ ranges over propositional variables, $A,A_i,B$ range over
formulas, $\,\Gamma$ ranges over disjoint unions of formulas, and
comma denotes disjoint union.\footnote{We label the conjunction and
disjunction rules with $\with$ and $\plus$ for reasons which will
become apparent later.}
By defining a sequent as a disjoint union of formulas, rather than an
ordered list, we avoid an exchange/permutation rule (\cf\
\cite[\S1.1]{TS96}).  Negation is primitive on propositional variables
$P$, and extends to compound formulas by de Morgan
duality.\footnote{$\:\neg (A\vee B)\,=\,(\neg A)\wedge (\neg B)\:$ and
$\:\neg(A\wedge B)\,=\,(\neg A)\vee(\neg B)\:$.}

The structural rules, weakening $\weakening$ and contraction
$\contraction$, are absorbed in the following variant, a right-sided
formulation of the propositional part of the calculus of \cite{Ket44},
called \gsthreep in \cite{TS96}.\footnote{We label the disjunction
rule as $\footnoteparr$ to distinguish it from the disjunction rule
$\plus$ of $\gsonep$.  The notation is derived from linear logic
\cite{Gir87}.}
\begin{labelledbox}{System \gsthreep}{\begin{center}
\begin{prooftree}\thickness=.08em
\rule{0pt}{1.3ex}
\justifies \;\Gamma,P,\neg P\; \using \axlabel
\end{prooftree}
\hspace{8ex}
\begin{prooftree}\thickness=.08em
\Gamma, A\hspace{8ex}\Gamma, B
\justifies \;\Gamma, A\wedge B\; \using \withlabel
\end{prooftree}
\hspace{8ex}
\begin{prooftree}\thickness=.08em
\Gamma, A,B
\justifies \;\Gamma, A\vee B\; \using \parlabel
\end{prooftree}
\end{center}}\end{labelledbox}
The new axiom $\Gamma,P,\neg P$ amounts to the original axiom $P,\neg
P$ followed immediately by weakenings.
This paper presents a propositional classical sequent calculus $\corep$
which is also free of structural rules:
\begin{labelledbox}{System \corep}{\begin{center}
\begin{prooftree}\thickness=.08em
\rule{0pt}{1.3ex}
\justifies \;P,\neg P\; \using \axlabel
\end{prooftree}
\hspace{4ex}
\begin{prooftree}\thickness=.08em
\Gamma,\Delta, A\hspace{4ex}\Gamma,\Sigma, B
\justifies \;\Gamma,\Delta,\Sigma, A\wedge B\; \using \andlabel
\end{prooftree}
\hspace{4ex}
\begin{prooftree}\thickness=.08em
\Gamma, A,B
\justifies \;\Gamma, A\vee B\; \using \parlabel
\end{prooftree}
\hspace{4ex}
\begin{prooftree}\thickness=.08em
\Gamma, A_i
\justifies \;\Gamma, A_1\vee A_2\; \using \plusilabel
\end{prooftree}
\end{center}}\end{labelledbox}
A distinguishing feature of \corep is the \emph{blended conjunction
rule}\footnote{By analogy with \gsthree and \gsthreep in \cite{TS96},
we reserve the symbol $\h$ for a full system with quantifiers, and use
$\corep$ to denote the propositional system.  Following
\cite{TS96}, we treat cut separately.  To maximise emphasis on
the blended conjunction rule, we omit quantifiers and cut in this
paper.}

\begin{displaymath}
\begin{prooftree}\thickness=.08em
\Gamma,\Delta, A\hspace{4ex}\Gamma,\Sigma, B
\justifies \;\Gamma,\Delta,\Sigma, A\wedge B\; \using \andlabel
\end{prooftree}\end{displaymath}
which combines the standard context-sharing and context-splitting
conjunction rules:
\begin{center}\label{usual}
\begin{prooftree}\thickness=.08em
\Gamma, A\hspace{4ex}\Gamma, B
\justifies \;\Gamma, A\wedge B\; \using \withlabel
\end{prooftree}
\hspace{8ex}
\begin{prooftree}\thickness=.08em
\Delta, A\hspace{4ex}\Sigma, B
\justifies \;\Delta,\Sigma, A\wedge B\; \using \tensorlabel
\end{prooftree}\end{center}
We refer to $\corep$ as (cut-free propositional)
\emph{minimal sequent calculus}.
In contrast to \gsthreep, contraction and weakening are absorbed
symmetrically: both are absorbed into the conjunction rule, leaving
the axiom rule intact.

$\Cp$ is evidently sound, since each of its rules can be derived
(encoded) in $\gsonep$.  Theorem~\ref{completeness}
(page~\pageref{completeness}) is completeness for formulas:
a formula is valid iff it is derivable in $\Cp$.\footnote{Completeness
here refers specifically to formulas, not to sequents. 
Section~\ref{degrees} discusses completeness for sequents.}

\subsection{Minimality}

The blended conjunction rule $\wedge$
is critical for the liberation from structural rules:
Proposition~\ref{cpminus-incomplete} (page~\pageref{cpminus-incomplete})
shows that relaxing it to the union of the the two standard
conjunction rules $\with$ and $\tensor$ breaks
completeness.\footnote{In other words, if we remove the $\wedge$ rule
and add both the $\with$ and the $\tensor$ rules, the resulting system
fails to be complete.  The formula $\:\big((P\wedge Q)\vee(\dual
Q\wedge P)\big)\vee \dual P\:$ becomes underivable (see the proof of
Proposition~\ref{cpminus-incomplete}, page~\pageref{cpminus-incomplete}).}
The main theorem of the paper
(page~\pageref{minimality-theorem}) formalises the sense in which
$\Cp$ is a minimal complete core of classical sequent calculus:
\vspace{1ex}
\begin{labelledbox}
{Theorem~\ref{minimality-theorem}: Minimality}
\sl \hspace{2ex}A standard sequent calculus $S$ is complete {\it iff}\/
$S\contains\corep$.\\
\end{labelledbox}
Here $S\contains T$ (``$S$ contains $T$'') iff every rule of $T$ is
derivable in $S$, and a \emph{standard sequent calculus} is any
propositional sequent calculus with the axiom $\overline{P,\neg P}$
and any subset of the following \emph{standard rules}:
\begin{center}
\begin{prooftree}\thickness=.08em
\Gamma,A\hspace{4ex}\Gamma, B
\justifies \;\Gamma, A\wedge B\; \using \withlabel
\end{prooftree}
\hspace{9ex}
\begin{prooftree}\thickness=.08em
\Gamma, A,B
\justifies \;\Gamma, A\vee B\; \using \parlabel
\end{prooftree}
\hspace{11ex}
\begin{prooftree}\thickness=.08em
\Gamma
\justifies \;\Gamma, A \using \weakeninglabel
\end{prooftree}

\vspace*{4ex}
\begin{prooftree}\thickness=.08em
\Delta,A\hspace{4ex}\Sigma, B
\justifies \;\Delta,\Sigma, A\wedge B\; \using \tensorlabel
\end{prooftree}
\hspace{8ex}
\begin{prooftree}\thickness=.08em
\Gamma, A_i
\justifies \;\Gamma, A_1\vee A_2\; \using \plusilabel
\end{prooftree}
\hspace{8ex}
\begin{prooftree}\thickness=.08em
\Gamma, A,A
\justifies \;\Gamma, A\; \using \contractionlabel
\end{prooftree}\vspace*{1ex}
\end{center}

\section{Notation and terminology}

Formulas are built from literals (propositional variables
$P,Q,R\ldots$ and their formal complements $\dual P, \dual Q,\dual
R,\ldots$) by the binary connectives \defn{and} $\wedge$ and \defn{or}
$\vee$.
Define \defn{negation} or \defn{not} $\neg$ as an operation on
formulas (rather than as a connective): $\neg P=\dual{P}$ and $\neg
\dual{P}=P$ for all propositional variables $P$, 
with $\neg (A\wedge B)=(\neg A)\vee (\neg B)$ and $\neg(A\vee B)=(\neg
A)\wedge(\neg B)$.

We identify a formula with its parse tree, a tree labelled with
literals at the leaves and connectives at the internal vertices.  
A \defn{sequent} is a non-empty disjoint union of
formulas.\footnote{Thus a sequent is a particular kind of labelled
forest.  This foundational treatment of formulas and sequents as
labelled trees and forests sidesteps the common problem of
``formulas'' versus ``formula occurrences'': disjoint unions of graphs
are well understood in graph theory \cite{Bol02}.}  Comma denotes
disjoint union.
Throughout the document, $P,Q,\ldots$ range over propositional
variables, $A,B,\ldots$ over formulas, and $\Gamma,\Delta,\ldots$ over
(possibly empty) disjoint unions of formulas.  

A formula $A$ is \defn{valid} if it evaluates to $1$ under all
possible $0/1$-assignments of its propositional variables (with the
usual interpretation of $\wedge$ and $\vee$ on $\{0,1\}$).
A sequent $A_1,\ldots,A_n$ is valid iff the formula $A_1\vee(A_2\vee (\ldots\vee (A_{n-1}\vee
A_n)\ldots))$ is valid.
A \defn{subsequent} of a sequent $\Gamma$ is any result of deleting
zero or more formulas from $\Gamma$; if at least one formula is
deleted, the result is a \defn{proper subsequent}.

\section{Completeness}

\begin{theorem}[Completeness]\label{completeness}
Every valid formula is derivable in $\Cp$.
\end{theorem}
The proof is via the following auxiliary definitions and lemmas.

A sequent is \defn{minimally valid}, or simply \defn{minimal}, if it
is valid while no proper subsequent is valid.  For example, the
sequents $P,\neg P$ and $\:P\wedge Q,\:\dual Q\wedge P,\:\dual P$ are
minimal, while $P,\neg P,Q$ is not.
\begin{lemma}\label{has-minvalid}
Every valid sequent contains a minimal subsequent.
\end{lemma}
\begin{proof}
Immediate from the definition of minimality.
\end{proof}
\begin{lemma}\label{literals}
Suppose a sequent $\:\Gamma\,$ is a disjoint union of literals (\ie,
$\:\Gamma\,$ contains no $\wedge$ or $\vee$).
Then $\;\Gamma\:$ is minimal iff $\;\Gamma\,=\,P,\,\neg P\:$ for some
propositional variable $P$.
\end{lemma}
\begin{proof}
By definition of validity in terms of valuations, $\:\Gamma\:$ is valid
iff it contains a complementary pair of literals, 
\ie, iff $\;\Gamma\,=\,P,\,\neg P,\,\Delta\;$ with $\,\Delta\,$ a disjoint union of
zero or more literals.  
Since $\,P,\,\neg P\,$ is valid, $\:\Gamma\:$ is minimal iff
$\,\Delta\,$ is empty.
\end{proof}
Suppose $\Gamma$ and $\Delta$ are each disjoint unions of formulas (so each
is either a sequent or empty).  Write $\Gamma\subforest\Delta$ if
$\Gamma$ results from deleting zero or more formulas from $\Delta$.
\begin{lemma}\label{and-union}
Suppose $\:\Gamma,\:A_1\wedge A_2\:$ is minimal.
Choose $\:\Gamma_1\subseteq\Gamma\:$ and $\:\Gamma_2\subseteq\Gamma\:$
such that $\;\Gamma_1,A_1$ and $\;\Gamma_2,A_2\:$ are minimal
(existing by Lemma~\ref{has-minvalid}, since $\;\Gamma,A_1\:$ and
$\;\Gamma,A_2\:$ are valid).  Then
every formula of $\:\Gamma\,$ is in at least one of the
$\:\Gamma_i\,$.
\end{lemma}
\begin{proof}
Suppose the formula $B$ of $\Gamma$ is in neither $\Gamma_i$.  Let
$\Gamma'$ be the result of deleting $B$ from $\Gamma$.
Then $\Gamma',\,A_1\wedge A_2\,$ is a valid 
proper subsequent of $\,\Gamma,\:A_1\wedge A_2\,$, contradicting
minimality.  
(The sequent $\,\Gamma',\,A_1\wedge A_2\,$
is valid since $\,\Gamma_1,A_1\,$ and $\,\Gamma_2,A_2\,$ are valid.)
\end{proof}
\begin{lemma}\label{or-sub}
Suppose $\;\Gamma\mkern-2mu,\:A\vee B\:$ is minimal and
$\;\Gamma\mkern-3mu,A\:$ is valid.  Then $\;\Gamma\mkern-3mu,A\:$
is minimal.
\end{lemma}
\begin{proof}
If not, some proper subsequent $\,\Delta\,$ of $\;\Gamma\mkern-2mu,A\:$ is valid.  If
$\,\Delta\,$ does not contain $A$, then it is also a proper subsequent of
$\;\Gamma\!,\,A\vee B\:$, contradicting minimality.  Otherwise let
$\,\Delta'\,$ be the result of replacing $A$ in $\,\Delta\,$ by $\:A\vee B\:$.
Since $\;\Delta\:$ is valid, so also is $\,\Delta'\,$.  Thus $\,\Delta'\,$ is a
valid proper subsequent of $\;\Gamma\!,\,A\vee B\,$, contradicting minimality.
\end{proof}
\begin{lemma}\label{par-sub}
Suppose $\;\Gamma,\,A\vee B\:$ is minimal and neither
$\;\Gamma\!,A\:$ nor $\;\Gamma\!,B\:$ is valid.  Then
$\;\Gamma\!,A,B\:$ is minimal.
\end{lemma}
\begin{proof}
Suppose $\;\Gamma,A,B\:$ had a valid proper subsequent $\,\Delta\,$.  Since
neither $\;\Gamma,A\:$ nor $\;\Gamma,B\:$ is valid, $\,\Delta\,$ must contain both
$A$ and $B$.  Let $\,\Delta'\,$ result from replacing $\,A,B\,$ by $\:A\vee B\:$
in $\,\Delta\,$.  Then $\,\Delta'\,$ is a valid proper subsequent of
$\:\Gamma,\,A\vee B\,$, contradicting minimality.
\end{proof}

Since a formula (viewed as a singleton sequent) is a minimal sequent,
the Completeness Theorem (Theorem~\ref{completeness}) is a special
case of:
\begin{proposition}\label{prop-minvalid}
Every minimal sequent is derivable in\/ $\Cp$.
\end{proposition}
\begin{proof}
Suppose $\:\Gamma\,$ is a minimal sequent.  We proceed by induction
on the number of connectives in $\Gamma$.
\begin{itemize}
\item \emph{Induction base (no connective).}
Since $\Gamma$ is minimal, Lemma~\ref{literals} implies $\Gamma=P,\neg
P$, the conclusion of the axiom rule $\overline{P,\neg P}$.
\item \emph{Induction step (at least one connective).}
\begin{enumerate}
\item \emph{Case: $\Gamma=\Delta,A_1\wedge A_2$.}  
By Lemma~\ref{and-union}, $\Gamma=\Sigma,\Delta_1,\Delta_2,A_1\wedge
A_2$ for $\Sigma,\Delta_1,A_1$ and $\Sigma,\Delta_2,A_2$ minimal.  Write down the conjunction rule
$$\begin{prooftree}\thickness=.08em
\Sigma,\Delta_1,A_1\hspace{4ex}\Sigma,\Delta_2,A_2
\justifies \;\Sigma,\Delta_1,\Delta_2, A_1\wedge A_2\; \using \andlabel
\end{prooftree}$$
and appeal to induction with the two hypothesis sequents.
\item \emph{Case: $\Gamma=\Delta,A_1\vee A_2$.}  
\begin{enumerate}
\item \emph{Case: $\Delta,A_i$ is valid for some $i\in\{1,2\}$.}  Write down the disjunction rule
$$\begin{prooftree}\thickness=.08em
\Delta, A_i
\justifies \;\Delta, A_1\vee A_2\; \using \plusilabel
\end{prooftree}$$
then appeal to induction with $\Delta,A_i$, which is minimal
by Lemma~\ref{or-sub}.
\item \emph{Case: $\Delta,A_i$ is not valid for each $i\in\{1,2\}$.}  
Thus $\Delta,A_1,A_2$ is minimal, by Lemma~\ref{par-sub}.
Write down the disjunction rule
$$\begin{prooftree}\thickness=.08em
\Delta, A_1,A_2
\justifies \;\Delta, A_1\vee A_2\; \using \parlabel
\end{prooftree}$$
then appeal to induction with $\Delta,A_1,A_2$.
\end{enumerate}
\end{enumerate}
\end{itemize}
($\Gamma$ may match both 1 and 2 in the inductive step,
permitting some choice in the construction of the derivation.  There is
choice in case 2(a) if both $\Delta,A_1$ and $\Delta,A_2$ are valid.)
\end{proof}

Note that completeness
does not hold for arbitrary valid sequents.  For example, the sequent
$\;P,\neg P,Q$ is valid but not derivable in $\Cp$.
A sequent is valid iff some some subsequent is derivable in $\Cp$.
Thus $\Cp$ is complete for sequents modulo final weakenings.  In this
sense, $\Cp$ is akin to system $\gsfivep$ of
\cite[\S7.4]{TS96} (related to resolution).
(See also Section~\ref{degrees}.)

\section{The Minimality Theorem}

Relaxing blended conjunction to the pair of standard conjunction rules
(context-sharing $\with$ and context-splitting $\tensor$) breaks
completeness.
Let $\Cpminus$ be the following subsystem of $\Cp$:\footnote{This
precursor of \Cp is (cut-free) multiplicative-additive linear logic
\cite{Gir87} with \emph{tensor} $\tensor$ and
\emph{with} $\with$ collapsed to $\wedge$, and \emph{plus} $\plus$ and \emph{par} $\footnoteparr$
collapsed to $\vee$.}\vspace*{1ex}
\begin{labelledbox}{System \Cpminus}{\begin{center}\label{cpminus}
\begin{prooftree}\thickness=.08em
\justifies \;P,\neg P\; \using \axlabel
\end{prooftree}
\hspace{7ex}\begin{minipage}[r]{3in}\begin{center}
\begin{prooftree}\thickness=.08em
\Gamma,A\hspace{4ex}\Gamma, B
\justifies \;\Gamma, A\wedge B\; \using \withlabel
\end{prooftree}
\hspace{9ex}
\begin{prooftree}\thickness=.08em
\Gamma, A,B
\justifies \;\Gamma, A\vee B\; \using \parlabel
\end{prooftree}\hspace{2ex}

\vspace*{3ex}
\begin{prooftree}\thickness=.08em
\Delta,A\hspace{4ex}\Sigma, B
\justifies \;\Delta,\Sigma, A\wedge B\; \using \tensorlabel
\end{prooftree}
\hspace{8ex}
\begin{prooftree}\thickness=.08em
\Gamma, A_i
\justifies \;\Gamma, A_1\vee A_2\; \using \plusilabel
\end{prooftree}\end{center}\end{minipage}
\end{center}}\end{labelledbox}
\begin{proposition}\label{cpminus-incomplete}
System $\Cpminus$ is incomplete.
\end{proposition}
\begin{proof}
We show that the valid formula $\:A\,=\,\big((P\wedge Q)\vee(\dual
Q\wedge P)\big)\vee \dual P\:$ is not derivable in $\Cpminus$.
The placement of the two outermost $\vee$ connectives forces the last
two rules of a potential derivation to be disjunction rules.  Since
$\:P\wedge Q,\,\dual Q\wedge P,\,\dual P\:$ is minimal (no proper
subsequent is valid), the two disjunction rules must be $\parr$ rather
than $\plus$:
\begin{displaymath}
\begin{prooftree}\thickness=.08em
\[
 P\wedge Q,\,\dual Q\wedge P,\,\dual P
 \justifies
 (P\wedge Q)\vee(\dual Q\wedge P),\,\dual P
 \using \parlabel
\]
\justifies \;\big((P\wedge Q)\vee(\dual
Q\wedge P)\big)\vee \dual P\; \using \parlabel
\end{prooftree}
\end{displaymath}
It remains to show that $\:P\wedge Q,\,\dual Q\wedge P,\,\dual P\:$ is
not derivable in $\Cpminus$.
There are only two connectives, both $\wedge$, so the last rule must
be a conjunction.
\begin{enumerate}
\item \emph{Case: the last rule is a context-sharing $\with$-rule.}
\begin{enumerate}
\item \emph{Case: The last rule introduces $P\wedge Q$.}
\begin{displaymath}
\begin{prooftree}\thickness=.08em
P,\,\dual Q\wedge P,\dual P
\hspace{4ex}
Q,\,\dual Q\wedge P,\,\dual P
\justifies 
\; P\wedge Q,\,\dual Q\wedge P,\,\dual P
\; \using \withlabel
\end{prooftree}
\end{displaymath}
The left hypothesis $P,\,\dual Q\wedge P,\,\dual P$ cannot be derived
in $\Cpminus$, since there is no $Q$ to match the $\dual Q$ (and no weakening).
\item \emph{Case: The last rule introduces $\dual Q\wedge P$.}
The same as the previous case, by symmetry, and exchanging
$Q\leftrightarrow \dual Q$.
\end{enumerate}
\item \emph{Case: the last rule is a context-splitting $\tensor$-rule.}
\begin{enumerate}
\item \emph{Case: The last rule introduces $P\wedge Q$.}
\begin{displaymath}
\begin{prooftree}\thickness=.08em
P,\,\Gamma
\hspace{4ex}
Q,\,\Delta
\justifies 
\; P\wedge Q,\,\dual Q\wedge P,\,\dual P
\; \using \tensorlabel
\end{prooftree}
\end{displaymath}
We must allocate each of $\dual Q\wedge P$ and $\dual P$ either
to $\Gamma$ or to $\Delta$.  If $\dual Q\wedge P$ is in $\Gamma$, then
$\,P,\Gamma\,$ is not derivable in $\Cpminus$, since it contains no $Q$
to match the $\dual Q$.  So $\dual Q\wedge P$ is in $\Delta$.  But
then the $\dual P$ is required in both $\Gamma$ and $\Delta$.
\item \emph{Case: The last rule introduces $\dual Q\wedge P$.}
The same as the previous case, by symmetry, and exchanging
$Q\leftrightarrow \dual Q$.
\end{enumerate}
\end{enumerate}
\end{proof}

A \defn{standard system}\label{standard-system} is any
propositional sequent calculus containing the axiom $\overline{P,\neg
P}$ and any of the following \defn{standard rules}:
\begin{center}
\begin{prooftree}\thickness=.08em
\Gamma,A\hspace{4ex}\Gamma, B
\justifies \;\Gamma, A\wedge B\; \using \withlabel
\end{prooftree}
\hspace{9ex}
\begin{prooftree}\thickness=.08em
\Gamma, A,B
\justifies \;\Gamma, A\vee B\; \using \parlabel
\end{prooftree}
\hspace{11ex}
\begin{prooftree}\thickness=.08em
\Gamma
\justifies \;\Gamma, A \using \weakeninglabel
\end{prooftree}

\vspace*{2ex}
\begin{prooftree}\thickness=.08em
\Delta,A\hspace{4ex}\Sigma, B
\justifies \;\Delta,\Sigma, A\wedge B\; \using \tensorlabel
\end{prooftree}
\hspace{8ex}
\begin{prooftree}\thickness=.08em
\Gamma, A_i
\justifies \;\Gamma, A_1\vee A_2\; \using \plusilabel
\end{prooftree}
\hspace{8ex}
\begin{prooftree}\thickness=.08em
\Gamma, A,A
\justifies \;\Gamma, A\; \using \contractionlabel
\end{prooftree}
\end{center}
Thus there are $2^6=64$ such systems (many of which will not be
complete).  

System $S$ \defn{contains} system $T$, denoted $S\contains T$, if each
rule of $T$ is a derived rule of $S$.  For example, system $\gsonep$
(page~\pageref{gsonep}) contains $\Cp$ since the blended conjunction
rule $\wedge$ and the disjunction rule $\parr$ of $\Cp$ can be derived
in $\gsonep$:
\begin{displaymath}
\begin{array}{ccc}
\begin{prooftree}\thickness=.08em
\Gamma,\Delta, A\hspace{4ex}\Gamma,\Sigma, B
\justifies \;\Gamma,\Delta,\Sigma, A\wedge B\; \using \andlabel
\end{prooftree}
&
\hspace{5ex}
\longleftarrow\label{wedge-from-Wwith}
\hspace{5ex}
&
\begin{prooftree}\thickness=.08em
\[
  \Gamma,\Delta, A 
  \justifies
  \Gamma,\Delta,\Sigma,A
  \using \weakeninglabel^\ast
\]
\hspace{4ex}
\[
  \Gamma,\Sigma, B
  \justifies
  \Gamma,\Delta,\Sigma,B
  \using \weakeninglabel^\ast
\]
\justifies \;\Gamma,\Delta,\Sigma, A\wedge B\; \using \withlabel
\end{prooftree}
\\\\\\
\begin{prooftree}\thickness=.08em
\Gamma, A,B
\justifies \;\Gamma, A\vee B\; \using \parlabel
\end{prooftree}
&
\hspace{5ex}
\longleftarrow
\hspace{5ex}
&
\begin{prooftree}\thickness=.08em
\[
  \[
    \Gamma, A,B
    \justifies \;\Gamma,A,A\vee B
    \using \plustwolabel
  \]
  \justifies \Gamma,A\vee B,A\vee B
  \using \plusonelabel
\]
\justifies \;\Gamma, A\vee B\; \using \contractionlabel
\end{prooftree}
\end{array}
\end{displaymath}
where $\weakeninglabel^\ast$ denotes a sequence of zero or more
weakenings.
\begin{theorem}[Minimality Theorem]\label{minimality-theorem}
A standard system is complete {\it iff}\/ it contains\/
$\Cp$.
\end{theorem}

\subsection{Proof of the Minimality Theorem}

Two systems are \defn{equivalent} if each contains the other.
For example, it is well known that $\gsonep$ (page~\pageref{gsonep})
is equivalent to:\footnote{This system is multiplicative
linear logic \cite{Gir87} plus contraction and weakening (with the
connectives denoted $\wedge$ and $\vee$ instead of $\tensor$ and
$\parr$).}
\begin{displaymath}
\begin{prooftree}\thickness=.08em
\rule{0pt}{1.3ex}
\justifies \;P,\neg P\; \using \axlabel
\end{prooftree}
\hspace{5ex}
\begin{prooftree}\thickness=.08em
\Delta, A\hspace{3ex}\Sigma, B
\justifies \;\Delta,\,\Sigma,\, A\wedge B\; \using \tensorlabel
\end{prooftree}
\hspace{5ex}
\begin{prooftree}\thickness=.08em
\Gamma, A_1, A_2
\justifies \;\Gamma,\, A_1\vee A_2\; \using \parlabel
\end{prooftree}
\hspace{5ex}
\begin{prooftree}\thickness=.08em
\Gamma
\justifies \;\Gamma, A \using \weakeninglabel
\end{prooftree}
\hspace{5ex}
\begin{prooftree}\thickness=.08em
\Gamma, A,A
\justifies \;\Gamma, A\; \using \contractionlabel
\end{prooftree}
\end{displaymath}
via the following rule derivations:
\begin{displaymath}
\begin{array}{ccc}
\begin{prooftree}\thickness=.08em
\Delta, A\hspace{4ex}\Sigma, B
\justifies \;\Delta,\Sigma, A\wedge B\; \using \tensorlabel
\end{prooftree}
&
\hspace{5ex}
\longleftarrow
\hspace{5ex}
&
\begin{prooftree}\thickness=.08em
\[
  \Delta, A 
  \justifies
  \Delta,\Sigma,A
  \using \weakeninglabel^\ast
\]
\hspace{4ex}
\[
  \Sigma, B
  \justifies
  \Delta,\Sigma,B
  \using \weakeninglabel^\ast
\]
\justifies \;\Delta,\Sigma, A\wedge B\; \using \withlabel
\end{prooftree}
\\\\\\
\begin{prooftree}\thickness=.08em
\Gamma, A_i
\justifies \;\Gamma, A_1\vee A_2\; \using \plusilabel
\end{prooftree}
&
\hspace{5ex}
\longleftarrow
\hspace{5ex}
&
\begin{prooftree}\thickness=.08em
\[
  \Gamma, A_i
  \justifies \;\Gamma,A_1,A_2
  \using \weakeninglabel
\]
\justifies \;\Gamma, A_1\vee A_2\; \using \parlabel
\end{prooftree}
\\\\\\
\begin{prooftree}\thickness=.08em
\Gamma, A,B
\justifies \;\Gamma, A\vee B\; \using \parlabel
\end{prooftree}
&
\hspace{5ex}
\longleftarrow
\hspace{5ex}
&
\begin{prooftree}\thickness=.08em
\[
  \[
    \Gamma, A,B
    \justifies \;\Gamma,A,A\vee B
    \using \plustwolabel
  \]
  \justifies \Gamma,A\vee B,A\vee B
  \using \plusonelabel
\]
\justifies \;\Gamma, A\vee B\; \using \contractionlabel
\end{prooftree}
\\\\\\
\begin{prooftree}\thickness=.08em
\Gamma, A\hspace{4ex}\Gamma, B
\justifies \;\Gamma, A\wedge B\; \using \withlabel
\end{prooftree}
&
\hspace{5ex}
\longleftarrow
\hspace{5ex}
&
\begin{prooftree}\thickness=.08em
\[
  \Gamma,A\hspace{4ex}\Gamma,B
  \justifies
  \;\Gamma,\Gamma, A\wedge B\; \using \tensorlabel
\]
\justifies \;\Gamma, A\wedge B\; \using \contractionlabel^\ast
\end{prooftree}
\end{array}
\end{displaymath}
We shall abbreviate these four rule derivations as follows, and write
analogous abbreviations for other rule
derivations.\footnote{``Context-splitting conjunction $\tensor$ is
derivable from context-sharing conjunction $\with$ and weakening
$\weakening$'', \etc}\vspace*{.7ex}
\begin{displaymath}\label{standard-derived-rules}
\begin{array}{ccc@{\hspace{12ex}}ccc}
\tensor & \longleftarrow & \with\,\weakening
&
\with   & \longleftarrow & \tensor\:\contraction
\\
\plus   & \longleftarrow & \parr\,\weakening
&
\parr   & \longleftarrow & \plus\:\contraction
\end{array}
\end{displaymath}

\subsubsection{The three complete standard systems}

As a stepping stone towards the Minimality Theorem, we shall prove
that, up to equivalence, there are only three complete standard
systems.

We abbreviate a system by listing its non-axiom rules.  For example,
$\gsonep=(\with,\plus,\weakening,\contraction)$ and
$\Cp=(\wedge,\plus,\parr)$.
Besides $\gsonep$, we shall pay particular attention to the systems
\begin{displaymath}
\begin{array}{r@{\;\;=\;\;}l@{\hspace{3ex}}l}
\positive & (\tensor,\plus,\contraction) & \text{\defn{Positive calculus}}\\[1ex]
\negative & (\with,\parr,\weakening) & \text{\defn{Negative calculus}}
\end{array}
\end{displaymath}
(Our terminology comes from polarity of connectives in linear
logic \cite{Gir87}: \emph{tensor} $\tensor$ and
\emph{plus} $\plus$ are positive, and \emph{with} $\with$ and
\emph{par} $\footnoteparr$ are negative.)
\begin{proposition}\label{only-ones}
Up to equivalence:
\begin{itemize}
\myitem{1}
\textnormal{$\gsonep\,=\,(\with,\plus,\contraction,\weakening)$}
is the only complete standard system with both contraction\/
$\contraction$ and weakening\/ $\weakening$;
\myitem{2} 
$\positive\,=\,(\tensor,\plus,\contraction)$ is the only complete
standard system without weakening\/ $\weakening$;
\myitem{3}  
$\negative\,=\,(\with,\parr,\weakening)$ is the only complete standard
system without contraction\/ $\contraction$.
\end{itemize}
\end{proposition}
The proof is via the following lemmas.
\begin{lemma}\label{containments}
$\Cp=(\wedge,\plus,\parr)$ is
contained in each of\/ $\positive\,=\,(\tensor,\plus,\contraction)$,
$\negative\,=\,(\with,\parr,\weakening)$ and\/
$\gsonep\,=\,(\with,\plus,\contraction,\weakening)$.
\end{lemma}
\begin{proof}
$\positive$ contains $\Cp$ since 
$\wedge\longleftarrow\contraction\tensor$\;,
\begin{displaymath}\label{wedge-from-Ctensor}
\begin{array}{ccc}
\begin{prooftree}\thickness=.08em
\Gamma,\Delta, A\hspace{4ex}\Gamma,\Sigma, B
\justifies \;\Gamma,\Delta,\Sigma, A\wedge B\; \using \andlabel
\end{prooftree}
&
\hspace{5ex}
\longleftarrow
\hspace{5ex}
&
\begin{prooftree}\thickness=.08em
\[
  \Gamma,\Delta,A
  \hspace{4ex}
  \Gamma,\Sigma,B
  \justifies
  \;\Gamma,\Gamma,\Delta,\Sigma, A\wedge B\;
  \using \tensorlabel
\]
\justifies \;\Gamma,\Delta,\Sigma, A\wedge B\; \using \contractionlabel^\ast
\end{prooftree}
\end{array}
\end{displaymath}
(where $\contraction^\ast$ denotes zero or more consecutive
contractions) and 
$\parr\longleftarrow\plus\contraction$\;:
\begin{displaymath}\label{derive}
\begin{array}{ccc}
\begin{prooftree}\thickness=.08em
\Gamma, A,B
\justifies \;\Gamma, A\vee B\; \using \parlabel
\end{prooftree}
&
\hspace{5ex}
\longleftarrow
\hspace{5ex}
&
\begin{prooftree}\thickness=.08em
\[
  \[
    \Gamma, A,B
    \justifies \;\Gamma,A,A\vee B
    \using \plustwolabel
  \]
  \justifies \Gamma,A\vee B,A\vee B
  \using \plusonelabel
\]
\justifies \;\Gamma, A\vee B\; \using \contractionlabel
\end{prooftree}
\end{array}
\end{displaymath}
$\negative$ contains $\Cp$ since 
$\wedge\longleftarrow\weakening\with$\;,
and
$\plus\longleftarrow\weakening\parr$
 (see page~\pageref{standard-derived-rules}).
$\gsonep\,=\,(\with,\plus,\contraction,\weakening)$ is equivalent to
$(\tensor,\with,\plus,\parr,\contraction,\weakening)$ since $\tensor$
and $\parr$ are derivable.  Thus $\gsonep$ contains $\positive$ (and
$\negative$), hence $\Cp$.
\end{proof}
\begin{lemma}\label{two-complete}
$\positive\,=\,(\contraction,\tensor,\plus)$ and\/
$\negative\,=\,(\with,\parr,\weakening)$ are complete.
\footnote{Recall that completeness refers to formulas, not sequents in general.}
\end{lemma}
\begin{proof}
Each contains $\Cp$ by Lemma~\ref{containments}, which is complete
(Theorem~\ref{completeness}).
\end{proof}
\begin{lemma}\label{one-CW}
Up to equivalence,
system
$\gsonep\,=\,(\with,\plus,\contraction,\weakening)$ is the only complete
standard system with both contraction $\contraction$ and weakening
$\weakening$.
\end{lemma}
\begin{proof}
$\gsonep$ is complete (see \eg\ \cite{TS96}, or by the fact that
$\gsonep$ contains $\Cp$ which is complete).  Any complete system must
have a conjunction rule ($\tensor$ or $\with$) and a disjunction rule
($\plus$ or $\parr$).  In the presence of $\contraction$ and
$\weakening$, the two conjunctions are derivable from one other, as
are the two disjunctions (see page~\pageref{standard-derived-rules}).
\end{proof}
\begin{lemma}\label{without-W}
A complete standard system without weakening\/ $\weakening$ must
contain\/
$\positive\,=\,(\tensor,\plus,\contraction)$.
\end{lemma}
\begin{proof}
System $\Cpminus\:=\:(\tensor,\plus,\with,\parr)$, with both
conjunction rules and both disjunction rules, is incomplete
(Proposition~\ref{cpminus-incomplete},
page~\pageref{cpminus-incomplete}), therefore we must have contraction
$\contraction$.

Without the $\plus$ rule, the valid formula $(P\vee \dual P)\vee Q$ is
not derivable: the last rule must be $\parr$, leaving us to derive
$P\vee\dual P,Q$, which is impossible without weakening $\weakening$
(\ie, with at most $\parr$, $\with$, $\tensor$ and $\contraction$
available), since, after a necessary axiom $\overline{P,\dual P}$ at
the top of the derivation, there is no way to introduce the formula $Q$.

Without the context-splitting $\tensor$ rule, the valid formula $P\vee
(Q\vee (\dual P\wedge \dual Q))$ is not derivable.  The last two rules
must be $\parr$, for if we use a $\plus$ we will not be able to match
complementary literals in the axioms at the top of the derivation.  Thus we
are left to derive $P,Q,\dual P\wedge \dual Q$, using $\with$ and
$\contraction$.
The derivation must contain an axiom rule $\overline{P,\dual P}$.  The next
rule can only be a $\with$ (since $P,\dual P$ cannot be the hypothesis
sequent of a contraction $\contraction$ rule).  Since the only
$\wedge$-formula in the final concluding sequent $P,Q,\dual P\wedge
\dual Q$ is $\dual P\wedge\dual Q$, and the $\with$ rule is context
sharing, the $\with$-rule must be
\begin{displaymath}
\begin{prooftree}\thickness=.08em
\[
  \justifies
  \;P,\dual P\;
\]
\hspace{4ex}
\[
  \proofdotseparation=1.2ex                 
  \proofdotnumber=4                         
  \leadsto                                  
  \;P,\dual Q
\]
\justifies \;P,\dual P\wedge \dual Q \; \using \withlabel
\end{prooftree}
\end{displaymath}
but $P,\dual Q$ is not derivable.
\end{proof}
\begin{lemma}\label{one-W}
Up to equivalence,
\textnormal{$\positive\,=\,(\tensor,\plus,\contraction)$} is the only complete standard system without
weakening\/ $\weakening$.
\end{lemma}
\begin{proof}
By Lemma~\ref{two-complete}, $\positive$ is complete.
By Lemma~\ref{without-W}, every $\weakening$-free complete standard
system contains $\positive$.
All other $\weakening$-free standard systems containing $\positive$
are equivalent to $\positive$, since the standard rule derivations
$\with\leftarrow\tensor\contraction$ and
$\parr\leftarrow\plus\contraction$ yield $\with$ and $\parr$ (see
page~\pageref{standard-derived-rules}).
\end{proof}
\begin{lemma}\label{without-C}
\!A complete standard system without contraction\/ $\contraction$ must
contain\/
$\negative=(\with,\!\parr,\!\weakening)$.
\end{lemma}
\begin{proof}
System $\Cpminus\:=\:(\tensor,\plus,\with,\parr)$, with both
conjunction rules and both disjunction rules, is incomplete
(Proposition~\ref{cpminus-incomplete},
page~\pageref{cpminus-incomplete}), therefore we must have weakening
$\weakening$.

Without the $\parr$ rule, the valid formula $P\vee \dual P$ would not be derivable.

Without the $\with$ rule 
the valid formula $P\vee (\dual P\wedge \dual P)$ would not be
derivable.
The last rule must be a $\parr$ (rather than a $\plus$, otherwise we
lack either $P$ or $\dual P$),
so we are left to derive $P,\,\dual P\wedge\dual P$.  The last rule
cannot be $\parr$ or $\plus$, as the only connective is $\wedge$.  It
cannot be $\weakening$, or else we lack either $P$ or $\dual P$.
It cannot be $\tensor$, as one of the two
hypotheses will be the single formula $\dual P$.
\end{proof}
\begin{lemma}\label{one-C}
Up to equivalence, 
$\negative\,=\,(\with,\parr,\weakening)$ is the only complete standard
system without contraction $\contraction$.
\end{lemma}
\begin{proof}
By Lemma~\ref{two-complete}, $\negative$ is complete.
By Lemma~\ref{without-C}, every $\contraction$-free complete standard
system contains $\negative$.
All other $\contraction$-free standard systems containing $\negative$
are equivalent to $\negative$, since the standard rule derivations
$\tensor\leftarrow\with\weakening$ and
$\plus\leftarrow\parr\weakening$ yield $\tensor$ and $\plus$ (see
page~\pageref{standard-derived-rules}).
\end{proof}
\begin{proofof}{Proposition~\ref{only-ones}}
Parts (1), (2) and (3) 
are Lemmas~\ref{one-CW}, \ref{one-W} and \ref{one-C}, respectively.
\end{proofof}
\begin{lemma}\label{no-WC-incomplete}
Every standard complete system has contraction\/ $\contraction$ or weakening\/
$\weakening$.
\end{lemma}

\begin{proof}
Otherwise it is contained in
$\Cpminus\,=\,(\tensor,\plus,\with,\parr)$, which is incomplete
(Prop.~\ref{cpminus-incomplete}).
\end{proof}
\begin{theorem}\label{only3}
Up to equivalence, there are only three complete standard systems: 
\begin{enumerate}
\item The
Gentzen-Sch\"utte-Tait system\/
$\gsonep\,=\,(\with,\plus,\contraction,\weakening)$.
\item
Positive calculus\/ $\positive\,=\,(\tensor,\plus,\contraction)$.
\item 
Negative calculus\/ $\negative\,=\,(\with,\parr,\weakening)$.
\end{enumerate}
\end{theorem}
\begin{proof}
Proposition~\ref{only-ones} and Lemma~\ref{no-WC-incomplete}.
\end{proof}

\begin{proofof}{Minimality Theorem (Theorem~\ref{minimality-theorem})}
Each of the three complete standard systems contains $\Cp$
(Lemma~\ref{containments}).
\end{proofof}

The three inequivalent complete standard systems $\gsonep$,
$\positive$ and $\negative$, together with propositional minimal sequent calculus
$\Cp$, sit in the following Hasse diagram of containments:
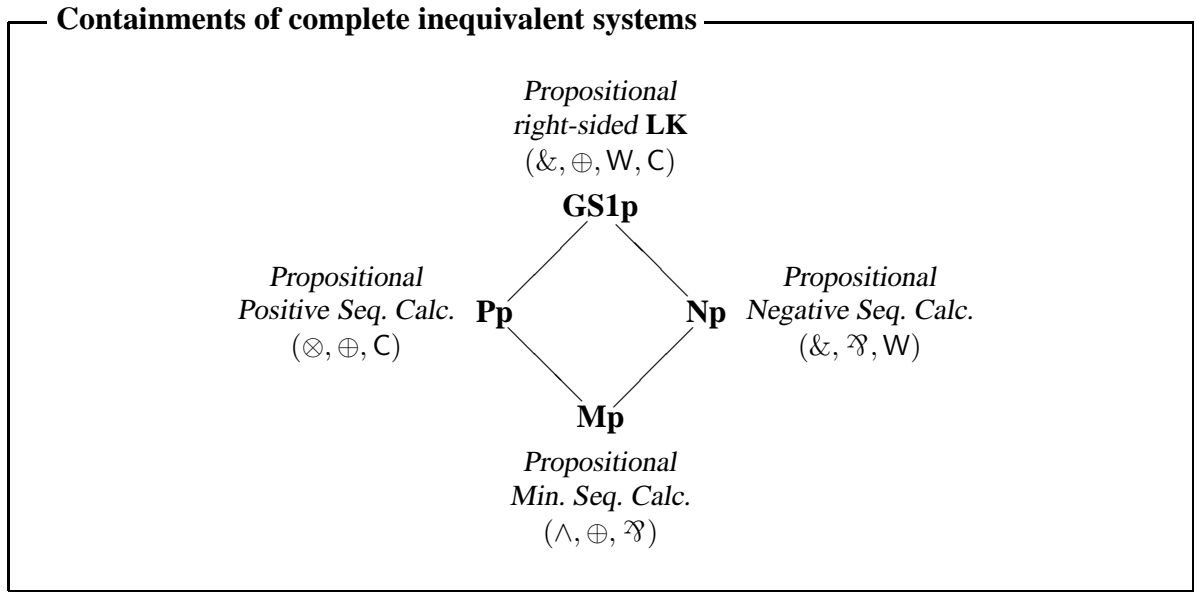
\begin{labelledbox}{Containments of complete inequivalent systems}
\begin{center}
\begin{picture}(0,185)(0,-90)
\put(-40,0){\makebox(0,0){$\positive$}\put(-15,0){\makebox(0,0)[r]
{\small\sl\shortstack{Propositional\\Positive Seq.\ Calc.\\$(\tensor,\plus,\contraction)$}}}}
\put(40,0){\makebox(0,0){$\negative$}\put(15,0){\makebox(0,0)[l]
{\small\sl\shortstack{Propositional\\Negative Seq.\ Calc.\\$(\with,\smallparr,\weakening)$}}}}
\put(0,40){\makebox(0,0){$\gsonep$}\put(0,12){\makebox(0,0)[b]
{\small\sl\shortstack{Propositional\\right-sided $\LK$\\$(\with,\plus,\weakening,\contraction)$}}}}
\put(0,-40){\makebox(0,0){$\Cp$}\put(0,-12){\makebox(0,0)[t]
{\small\sl\shortstack{Propositional\\Min.\ Seq.\ Calc.\\$(\wedge,\plus,\smallparr)$}}}}
\put(-35,5){\line(1,1){30}}
\put(33,7){\line(-1,1){27}}
\put(35,-5){\line(-1,-1){30}}
\put(-35,-5){\line(1,-1){28}}
\end{picture}
\end{center}
\end{labelledbox}
Thus we can view propositional minimal sequent calculus $\Cp$ as a minimal
complete core of $\gsonep$, hence of (propositional) Gentzen's
$\LK$.

\section{Extended Minimality Theorem}\label{ext-min}

Define an \defn{extended system} as one containing the axiom rule
$\overline{P,\dual P}$ and any of the following rules.  (We have
extended the definition of \emph{standard system} by making blended
conjunction available.)
\begin{labelledbox}{Extended system rules}\label{extended-rules}\begin{displaymath}
\begin{array}{c}
\\[-1ex]
\begin{prooftree}\thickness=.08em
\Delta,A\hspace{4ex}\Sigma, B
\justifies \;\Delta,\Sigma, A\wedge B\; \using \tensorlabel
\end{prooftree}
\hspace{6ex}
\begin{prooftree}\thickness=.08em
\Gamma,\Delta, A\hspace{4ex}\Gamma,\Sigma, B
\justifies \;\Gamma,\Delta,\Sigma, A\wedge B\; \using \andlabel
\end{prooftree}
\hspace{6ex}
\begin{prooftree}\thickness=.08em
\Gamma,A\hspace{4ex}\Gamma, B
\justifies \;\Gamma, A\wedge B\; \using \withlabel
\end{prooftree}
\\\\\\[-1ex]
\begin{prooftree}\thickness=.08em
\Gamma, A_i
\justifies \;\Gamma, A_1\vee A_2\; \using \plusilabel
\end{prooftree}
\hspace{12ex}
\begin{prooftree}\thickness=.08em
\Gamma, A,B
\justifies \;\Gamma, A\vee B\; \using \parlabel
\end{prooftree}
\\\\\\[-1ex]
\begin{prooftree}\thickness=.08em
\Gamma, A,A
\justifies \;\Gamma, A\; \using \contractionlabel
\end{prooftree}
\hspace{14ex}
\begin{prooftree}\thickness=.08em
\Gamma
\justifies \;\Gamma, A \using \weakeninglabel
\end{prooftree}
\\\\
\end{array}
\end{displaymath}
\end{labelledbox}
The Minimality Theorem (Theorem~\ref{minimality-theorem},
page~\pageref{minimality-theorem}) extends as follows.
\begin{theorem}[Extended Minimality Theorem]\label{extended-minimality-theorem}
An extended system is complete iff it contains propositional minimal
sequent calculus\/ $\Cp$.
\end{theorem}
To prove this theorem, we require two lemmas.
\begin{lemma}\label{S-with}
Suppose\/ $S$ is a complete extended system with the blended
conjunction rule\/ $\wedge$, and with at least one of contraction\/
$\contraction$ or weakening\/ $\weakening$.
Then $S$ is equivalent to a standard system.
\end{lemma}
\begin{proof}
If $S$ has weakening $\weakening$, let $S'$ be the result of replacing
the blended conjunction rule $\wedge$ in $S$ by context-sharing
conjunction $\with$; otherwise $S$ has contraction, and let $S'$
result from replacing $\wedge$ by context-splitting $\tensor$.  Then
$S'$ is equivalent to $S$, since $\wedge\longleftarrow\tensor\contraction$
(page~\pageref{wedge-from-Ctensor}) and 
$\wedge\longleftarrow\with\weakening$
(page~\pageref{wedge-from-Wwith}).
\end{proof}
\begin{lemma}\label{S-without}
Suppose\/ $S$ is a complete extended system with neither contraction\/
$\contraction$ nor weakening\/ $\weakening$.
Then $S$ is equivalent to propositional minimal sequent calculus\/ $\Cp$.
\end{lemma}
\begin{proof}
Since $\Cpminus\,=\,(\tensor,\with,\plus,\parr)$ is incomplete
(Proposition~\ref{cpminus-incomplete},
page~\pageref{cpminus-incomplete}), $S$ must have the blended
conjunction rule $\wedge$ either directly or as a derived rule.  Since
$S$ is complete, it must have a disjunction rule, therefore it could
only fail to be equivalent to $\Cp\,=\,(\wedge,\plus,\parr)$ if (a) it
has $\plus$ and $\parr$ is not derivable, \ie, $S$ is equivalent to
$(\wedge,\plus)$, or (b) it has $\parr$ and $\plus$ is not derivable,
\ie, $S$ is equivalent to $(\wedge,\parr)$.
In case (a), the valid formula $P\vee\dual P$ would not be derivable,
and in case (b) the valid formula $(P\vee\dual P)\vee Q$ would not be
derivable, either way contradicting the completeness of $S$.
\end{proof}
\begin{proofof}{the Extended Minimality Theorem (Theorem~\ref{extended-minimality-theorem})}
Suppose $S$ is a complete extended system.
If $S$ has contraction $\contraction$ or weakening $\weakening$ then
it is equivalent to a standard system by Lemma~\ref{S-with}, hence
contains $\Cp$ by the original Minimality Theorem.
Otherwise $S$ is equivalent to $\Cp$ by Lemma~\ref{S-without}, hence
in particular contains $\Cp$.

Conversely, suppose $S$ is an extended system containing $\Cp$.  Then
$S$ is complete since $\Cp$ is complete.
\end{proofof}
We also have the following extension of Theorem~\ref{only3}
(page~\pageref{only3}), which stated that, up to equivalence, there are
only three complete standard systems, $\gsonep$, $\positive$ and
$\negative$.
\begin{theorem}\label{only4}
Up to equivalence, there are only four complete extended systems:
\begin{enumerate}
\item The
Gentzen-Sch\"utte-Tait system\/
$\gsonep\,=\,(\with,\plus,\contraction,\weakening)$.
\item Positive calculus\/ $\positive\,=\,(\tensor,\plus,\contraction)$.
\item Negative calculus\/ $\negative\,=\,(\with,\parr,\weakening)$.
\item Propositional minimal sequent calculus\/ $\Cp\,=\,(\wedge,\plus,\parr)$.
\end{enumerate}
\end{theorem}
\begin{proof}
Theorem~\ref{only3} together with Lemmas~\ref{S-with} and
\ref{S-without}.
\end{proof}

\section{Degrees of completeness}\label{degrees}

We defined a system as \emph{complete} if every valid formula
(singleton sequent) is derivable.  To avoid ambiguity with forthcoming
definitions, let us refer to this default notion of completeness as
\defn{formula-completeness}.
Define a system as
\defn{minimal-complete} if every minimal\footnote{Recall that a 
valid sequent is minimal if no proper subsequent is valid.} sequent is
derivable, and \defn{sequent-complete} if every valid sequent is
derivable.  (Thus \emph{sequent-complete} implies
\emph{minimal-complete} implies \emph{formula-complete}.)

For a minimal-complete system $S$, a sequent $\Gamma$ is valid iff a
subsequent of $\Gamma$ is derivable in $S$.  Thus a minimal-complete
system $S$ can be viewed as sequent-complete, modulo final weakenings.
(\Cf\ system $\gsfivep$ of \cite[\S7.4]{TS96} (related to resolution).)

\begin{proposition}\label{not-seq-complete}
$\positive\,=\,(\tensor,\plus,\contraction)$ and\/
$\Cp\,=\,(\wedge,\plus,\parr)$ are formula-complete
and minimal-complete, but not sequent-complete.
\end{proposition}
\begin{proof}
We have already proved that $\Cp$ (hence also $\positive$, by
containment) is minimal-complete (Proposition~\ref{prop-minvalid}).

We show that the valid (non-minimal) sequent $P,\dual P,Q$ is
not derivable in $\positive$ (hence also in $\Cp$).  A derivation must contain
an axiom rule $\overline{P,\dual P}$.  This cannot be followed by a
$\tensor$ or $\plus$ rule, otherwise we introduce a connective
$\wedge$ or $\vee$ which cannot subsequently be removed by any other
rule before the concluding sequent $P,\dual P,Q$.  Neither can it be
followed by contraction $\contraction$, since there is nothing to
contract.
\end{proof}
\begin{proposition}\label{seq-complete}
$\negative\,=\,(\with,\parr,\weakening)$ is formula-,
minimal- and sequent-complete.
\end{proposition}
\begin{proof}
$\negative$ is minimal-complete since it contains $\Cp$.
Suppose $\Gamma$ is a valid but not minimal sequent.  Choose a minimal
subsequent $\Delta$ of $\Gamma$ (see Lemma~\ref{has-minvalid},
page~\pageref{has-minvalid}).  By minimal-completeness, $\Delta$ has a
derivation.  Follow this with weakenings to obtain $\Gamma$.
\end{proof}
Below we have annotated our Hasse diagram with completeness
strengths.
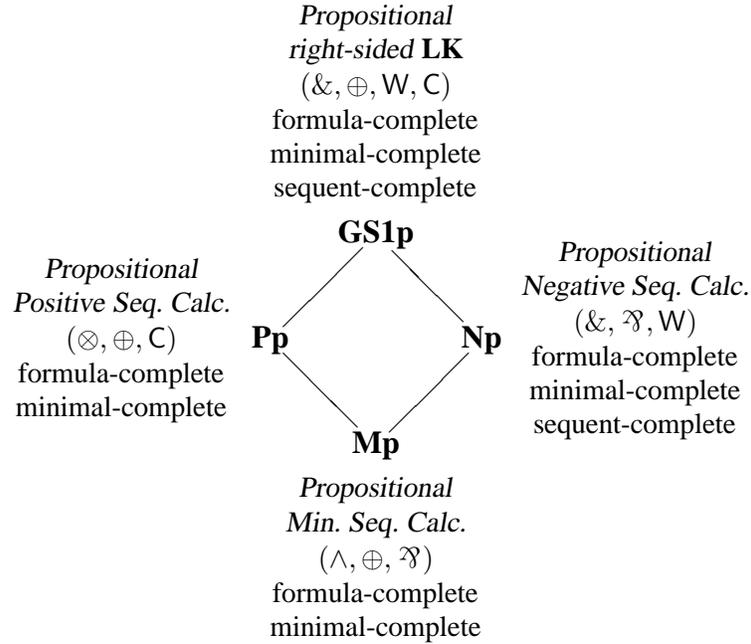
\vspace*{1ex}
\begin{center}
\begin{picture}(0,260)(0,-125)
\put(-40,0){\makebox(0,0){$\positive$}\put(-15,0){\makebox(0,0)[r]
{\small\sl\shortstack{Propositional\\Positive Seq.\ Calc.\\$(\tensor,\plus,\contraction)$
\\{}\normalfont formula-complete
\\{}\normalfont minimal-complete
}}}}
\put(40,0){\makebox(0,0){$\negative$}\put(15,0){\makebox(0,0)[l]
{\small\sl\shortstack{Propositional\\Negative Seq.\ Calc.\\$(\with,\smallparr,\weakening)$
\\{}\normalfont formula-complete
\\{}\normalfont minimal-complete
\\{}\normalfont sequent-complete
}}}}
\put(0,40){\makebox(0,0){$\gsonep$}\put(0,12){\makebox(0,0)[b]
{\small\sl\shortstack{Propositional\\right-sided $\LK$\\$(\with,\plus,\weakening,\contraction)$
\\{}\normalfont formula-complete
\\{}\normalfont minimal-complete
\\{}\normalfont sequent-complete
}}}}
\put(0,-40){\makebox(0,0){$\Cp$}\put(0,-12){\makebox(0,0)[t]
{\small\sl\shortstack{Propositional\\Min.\ Seq.\ Calc.\\$(\wedge,\plus,\smallparr)$
\\{}\normalfont formula-complete
\\{}\normalfont minimal-complete
}}}}
\put(-35,5){\line(1,1){30}}
\put(33,7){\line(-1,1){27}}
\put(35,-5){\line(-1,-1){30}}
\put(-35,-5){\line(1,-1){28}}
\end{picture}
\end{center}

\section{Possible future work}
\begin{enumerate}
\item \emph{Cut.}  
Chapter 4 of \cite{TS96} gives a detailed analysis of cut for Gentzen
systems.
One could pursue an analogous analysis of cut for minimal sequent
calculus.  Aside from context-splitting and context-sharing cut rules
\begin{center}
\begin{prooftree}\thickness=.08em
\Delta, A\hspace{4ex}\Sigma, \neg A
\justifies \;\Delta,\Sigma \; \using \splitcutlabel
\end{prooftree}
\hspace{8ex}
\begin{prooftree}\thickness=.08em
\Gamma, A\hspace{4ex}\Gamma, \neg A
\justifies \;\Gamma\; \using \supcutlabel
\end{prooftree}
\end{center}
one might also investigate a blended cut rule:
\begin{center}
\begin{prooftree}\thickness=.08em
\Gamma, \Delta, A\hspace{4ex}\Gamma, \Sigma, \neg A
\justifies \;\Gamma, \Delta, \Sigma\; \using \cutlabel
\end{prooftree}
\end{center}
\item \emph{Quantifiers.}  
Explore the various ways of adding quantifiers to $\Cp$, for a full first-order system $\h$.
\item \emph{Mix (nullary multicut).}  
Gentzen's multicut rule
\begin{center}
\begin{prooftree}\thickness=.08em
\Delta, A_1,\ldots,A_m\hspace{4ex}\Sigma, \neg A_1,\ldots, \neg A_n
\justifies \;\Delta,\Sigma \;
\end{prooftree}
\end{center}
in the nullary case $m=n=0$ has been of particular interest to linear
logicians \cite{Gir87}, who call it the \emph{mix} rule.  One could
investigate context-splitting, context-sharing and blended
incarnations:
\begin{displaymath}
\begin{prooftree}\thickness=.08em
\Delta\hspace{4ex}\Sigma
\justifies \;\Delta,\Sigma \; \using \splitmixlabel
\end{prooftree}
\hspace{8ex}
\begin{prooftree}\thickness=.08em
\Gamma\hspace{4ex}\Gamma
\justifies \;\Gamma\; \using \supmixlabel
\end{prooftree}
\hspace{8ex}
\begin{prooftree}\thickness=.08em
\Gamma, \Delta\hspace{4ex}\Gamma, \Sigma
\justifies \;\Gamma, \Delta, \Sigma\; \using \mixlabel
\end{prooftree}
\end{displaymath}
\end{enumerate}

\bibliographystyle{alpha}

\end{document}